\documentclass[12pt]{article}
\usepackage{amsfonts}
\usepackage{mathrsfs}
\usepackage{epsfig,latexsym,amsfonts,amsmath,amssymb}
\newcommand{\qed}{\hfill$\Box$}

\newtheorem{prob}[subsection]{Problem}
\newtheorem{proposition}[subsection]{Proposition}

\newtheorem{corollary}[subsection]{Corollary}

\newtheorem{lem}[subsection]{Lemma}
\newtheorem{thm}[subsection]{Theorem}

\newcommand{\vanish}[1]{}

\begin{document}
\parskip=12pt

\thispagestyle{empty} \setlength{\baselineskip}{1.2\baselineskip}

\title{Primitivity and Independent Sets in Direct Products of Vertex-Transitive Graphs
\thanks{Supported by the National
Natural Science Foundation of China (No.10826084) and Zhejiang
Innovation Project (Grant No. T200905).}}
\author{Huajun Zhang  \\
{\small Department of Mathematics, Shanghai Normal University,
Shanghai 200234, China}\\ {\small and} \\ {\small Department of
Mathematics,
 Zhejiang Normal University, Jinhua 321004, P.R. China }\\
  {\small E-mail: huajunzhang@zjnu.cn}}
\date{}

\maketitle

\noindent\textbf{Abstract.}  We introduce the concept of the
primitivity of independent set in vertex-transitive graphs, and
investigate the relationship between the primitivity and the
structure of maximum independent sets in direct products of
vertex-transitive graphs.  As a consequence of our main results, we
positively solve an open problem related to the structure of
independent sets in powers of vertex-transitive graphs.
\section{Introduction}

The direct product $G\times H$ of two graphs $G$ and $H$ is defined
by \[V(G\times H)=V(G)\times V(H)\] and \[E(G\times
H)=\{[(u_1,u_2),(v_1,v_2)]: [u_1,v_1]\in E(G) \mbox{\ and \ }
[u_2,v_2]\in E(H)\}.\] For a graph $G$, let $G^n=G\times \cdots
\times G$ denote the $n$-th power of $G$.

It is clear that if $I$ is an independent set of $G$ (or $H$), then
$I\times H$ (or $G\times I$) is an independent set of $G\times H$.
We say that  $G\times H$ is {\em MIS-normal}
(maximum-independent-set-normal) if each of its maximum independent
sets is of this form. Then the independence number
\begin{eqnarray}\label{eq12}\alpha(G\times
H)=\max\{\alpha(G)|H|, \alpha(H)|G|\}\end{eqnarray} if $G\times H$
is MIS-normal. A product $G_1\times G_2\times \cdots \times G_n$ is
said to be MIS-normal if all of its maximum independent sets are
preimages of projections of maximum independent sets of one of its
factors.

This poses two immediate problems: whether (1) holds for all graphs
$G$ and $H$, and whether $G\times H$ is MIS-normal when (1) holds.
 In general, however, (1) does not hold for some
non-vertex-transitive graphs (see \cite{Kla}). So, Tardif
\cite{Tardif} asked whether (1) holds for all vertex-transitive
graphs $G$ and $H$. Larose and Tardif \cite{larose} investigated the
relationship between the projectivity and the structure of maximal
independent sets in  powers of a circular graph, Kneser graph, or
truncated simplex. Recently, Mario and Vera \cite{MJ} proved that
(1) holds for some special vertex-transitive graphs, e.g., circular
graphs and Kneser graphs. In fact,  Frankl \cite{frankl3} proved in
1996, one year before Tardif's question was posed, that (1) holds
for Kneser graphs. Subsequently, Ahlswede, Aydinian and Khachatrian
\cite{Ah1} generalized Frankl's result.

In the context of vertex-transitive graphs, the ``No-Homomorphism"
lemma of Albertson and Collins \cite{makl} is useful to get bounds
on the size of independent sets.
\begin{lem}\label{ac}(\cite{makl})
Let $G$ and $H$ be two graphs such that $G$ is vertex-transitive and
there exists a homomorphism $\phi: H\mapsto G$. Then
$\frac{\alpha(G)}{|V(G)|}\leq \frac{\alpha(H)}{|V(H)|},$ and the
equality holds if and only if for any independent set $I$ of
cardinality $\alpha(G)$ in $G$, $\phi^{-1}(I)$ is an independent set
of  cardinality $\alpha(H)$ in $H$.
\end{lem}

By this lemma, it is easy to deduce that
$\alpha(G^n)=\alpha(G)|V(G)|^{n-1}$ for any vertex-transitive graph
$G$ and positive integer $n$ (see \cite{larose}). So it is natural
to ask whether $G^n$ is MIS-normal.   Evidently, if $G^n$ is
MIS-normal for some $n > 2$, so is $G^2$. Conversely, Larose and
Tardif \cite{larose} posed the following problem.
\begin{prob}\label{prb1}(see \cite{larose} )Let $G$ be a non-bipartite vertex-transitive graph.
If $G^2$ is MIS-normal, is the same for all powers of $G$?
\end{prob}

This paper is organized as follows. In the next section, we
introduce a concept  of the primitivity of independent sets in a
vertex-transitive graph, and prove that the primitivity can be
preserved in direct products under certain conditions. Based on
these results we  establish in section $3$ a direct product theorem
on the MIS-normality.
 As a consequence, Problem \ref{prb1} is
solved.
\section{ Primitivity of independent sets}
In the sequel of this paper, let $G$ and $H$ be vertex-transitive
graphs. By $I(G)$ we denote the set of all maximum independent sets
of $G$. For any subset $A$ of $V (G)$, let $\alpha(A)$ denote the
independence number of the induced subgraph of $G$ by $A$, and we
define
$$N_G(A)=\{b\in G:\mbox{$(a,b)\in E(G)$ for some $a\in A$}\},$$ $$
N_G[A]= N_G(A)\cup A\ \mbox{and}\ \overline{N}_G[A]=G-N_G[A].$$ In
Lemma \ref{ac}, by taking $H$  as an induced subgraph of $G$ and
$\phi$ as the embedding mapping, we obtain the following lemma (cf.
\cite{Cameron}).
\begin{lem}\label{caku}
$\frac{\alpha(G)}{|V(G)|}\leq \frac{\alpha(B)}{|B|}$ holds for all
$B\subseteq V(G)$. Equality implies that $|S\cap B|=\alpha(B)$ for
every $S\in I(G)$.
\end{lem}

A graph $G$ is said to be \textit{non-empty} if $E(G)\neq\emptyset$.
Lemma \ref{caku} implies that $\alpha(G)\leq |V(G)|/2$ for all
non-empty vertex-transitive graphs. Equality holds if and only if
$G$ is bipartite, which we state as a corollary  for reference.
\begin{corollary}\label{bi}
Let $G$ be a non-empty vertex-transitive graph. Then
$\frac{\alpha(G)}{|G|}\leq \frac{1}{2}$, and equality holds if and
only if $G$ is  bipartite.
\end{corollary}

\begin{proposition}\label{coro1}
Let $A$ be an independent set of $G$.  Then
$\frac{|A|}{|N_G[A]|}\leq \frac{\alpha(G)}{|V(G)|}$. Equality
implies that  $|S\cap N_G[A]|=|A|$  for every $S\in
  I(G)$, and in particularly $A\subseteq S$ for some $S\in I(G)$.
\end{proposition}
\textbf{Proof.} Since $A$ is an independent set, clearly
\[\frac{|A|+\alpha(\overline{N}_G[A])}{|N_G[A]|+|\overline{N}_G[A]|}\leq \frac{\alpha(G)}{|V(G)|}.\]
By Lemma \ref{caku} we see that
$\frac{\alpha(\overline{N}_G[A])}{|\overline{N}_G[A]|}\geq\frac{\alpha(G)}{|V(G)|}$,
so $\frac{|A|}{|N_G[A]|}\leq \frac{\alpha(G)}{|V(G)|}$. Equality in
the latter implies equality in the former. In this case any $S\in
I(G)$ must be the union of a maximum independent set in
$\overline{N}_G[A]$ and an independent set of size $|A|$ in
$N_G[A]$, and thus $|S\cap N_G[A]|=|A|$.\qed

An independent set $A$ in $G$ is said to be {\em imprimitive} if
$|A|<\alpha(G)$ and $\frac{|A|}{|N_G[A]|}=\frac{\alpha(G)}{|V(G)|}$.
We say that $G$ is {\em IS-imprimitive} if $G$ has an imprimitive
independent set. In the other case, $G$ is \emph{IS-primitive}.

\begin{proposition}\label{irregular}
Let $A$ be a maximum imprimitive independent set of $G$. Set
$B=\overline{N}_G[A]$.  Then
$\frac{\alpha(B)}{|B|}=\frac{\alpha(G)}{|V(G)|}$ and
  $\{\sigma(B)|\sigma\in\mbox{Aut}(G)\}$ forms a nontrivial partition of $V(G)$,
  i.e., $\sigma(B)\cap B=\emptyset$ or $B$ for each $\sigma\in \mbox{Aut}(G)$.
\end{proposition}
\textbf{Proof.} Clearly $\frac{|A|+\alpha(B)}{|N_G[A]|+|B|}\leq
\frac{\alpha(G)}{|V(G)|}$. Combining the condition of $A$ and Lemma
\ref{caku}, we have
$\frac{\alpha(B)}{|B|}=\frac{\alpha(G)}{|V(G)|}$. By definition,
$N_G[\sigma(A)]=\sigma(N_G[A])$ for all $\sigma\in \mbox{Aut}(G)$.
Suppose that there exists a $\sigma\in \mbox{Aut}(G)$ such that $
\sigma(B)\neq B$ and $\sigma(B)\cap B\neq \emptyset$. Then
$\sigma(N_G[A])\neq N_G[A]$ and
\begin{eqnarray}\label{ed}|V(G)|>|N_G[A]\cup \sigma\big(N_G[A]\big)|>|N_G[A]|.\end{eqnarray}
Let $C=\sigma(A)\cup (A-N_G[\sigma(A)])$. Then $C$ is also an
independent set and $$N_G[C]\subseteq N_G[A]\cup \sigma(N_G[A]).$$
By Proposition \ref{coro1},  $|S\cap N_G[A]|=|A|$ for all $S\in
I(G)$, which implies that $(S-N_G[A])\cup A\in I(G)$ for all $S\in
I(G)$. Similarly,
\begin{eqnarray*}
  &&((S-N_G[A])\cup A)-N_G[\sigma(A)])\cup \sigma(A)
\\
&=&(S-N_G[A]\cup N_G[\sigma(A)])\cup (A-N_G[\sigma(A)])\cup
\sigma(A)\\
&=&(S-N_G[A]\cup N_G[\sigma(A)])\cup C
\end{eqnarray*}
is also a maximum independent set of $G$, which implies $|S\cap
(N_G[A]\cup N_G[\sigma(A)])|=|C|$ for all $S\in I(G)$.

 Given a $u\in V(G)$, suppose that there are $r$ $S$'s in $I(G)$ such
that $u\in S$. Since $G$ is vertex-transitive,  the number $r$ is
independent of the choice of $u$. Thus $r|V(G)|=\alpha(G)|I(G)|$. On
the other hand, since  $|S\cap (N_G[A]\cup N_G[\sigma(A)])|=|C|$ for
all $S\in I(G)$, $|C||I(G)|=r|N_G[A]\cup N_G[\sigma(A)]|$. Combining
the above two equalities, we have $\frac{|C|}{|N_G[A]\cup
N_G[\sigma(A)]|}=\frac{\alpha(G)}{|V(G)|}$.
 Thus, by Proposition \ref{coro1} we have
$$\frac{\alpha(G)}{|V(G)|}\geq \frac{|C|}{N_G[C]}\geq
\frac{|C|}{|N_G[A]\cup N_G[\sigma(A)]|}=\frac{\alpha(G)}{|V(G)|},$$
which implies $N_G[C]=N_G[A]\cup N_G[\sigma(A)]$ and
$\frac{|C|}{|N_G[C]|}=\frac{\alpha(G)}{|V(G)|}$. By (\ref{ed}), we
have $|A|<|C|<\alpha(G)$, contradicting the maximality of $|A|$.
This completes the proof. \qed

 The concept of primitivity comes from permutation groups: A
permutation group $\Gamma$ acting on a set $X$ is called primitive
if $\Gamma$ preserves no nontrivial partition of $X$. In the other
case, $\Gamma$ is imprimitive. As usual (see e.g. \cite{larose}), a
vertex-transitive graph $G$ is called primitive if its automorphism
group, as a permutation group on $V(G)$, is primitive. By
Proposition \ref{irregular} we see that if $G$ is primitive, then
$G$ is IS-primitive. But the converse is not true.

For any
  $S\subseteq V(G)\times V(H)$, $a\in G$ and $u\in H$, define
  $$\partial_G(u,S)=\{b\in G: (b,u)\in S\}, \ \ \partial_H(a,S)=\{v\in H: (a,v)\in S\},$$
and
$$\partial_G(S)=\{b\in G: \partial_H(b,S)\neq\emptyset\}, \ \ \partial_H(S)=\{v\in H: \partial_G(v,S)\neq\emptyset\}.$$
By definition we see that $\partial_G(S)$ and $\partial_H(S)$ are in
fact the projections of $S$ on $G$ and $H$, respectively.

  \begin{lem}\label{dgp}
  Suppose $G\times H$ is MIS-normal and $\frac{\alpha(H)}{|H|}\leq \frac{\alpha(G)}{|G|}$.
   If $G\times H$ is  IS-imprimitive, then one of the following two possible cases holds:
\begin{itemize}
  \item [\rm(i)] $\frac{\alpha(H)}{|H|}=\frac{\alpha(G)}{|G|}$,  and  one of them is IS-imprimitive or both $G$ and $H$ are bipartite;
  \item [\rm(ii)]$\frac{\alpha(H)}{|H|}<\frac{\alpha(G)}{|G|}$, and   $G$ is IS-imprimitive or $H$ is disconnected.
\end{itemize}
\end{lem}
  \textbf{Proof.} Throughout this proof, we denote $N_{G\times H}[A]$
by $N[A]$ for brevity. Suppose that $G\times H$ is IS-imprimitive
and let $A$ be a maximum imprimitive independent set of $G\times H$.
Clearly, $\alpha(G\times H)=\alpha(G)|V(H)|$, and thus
$\frac{|A|}{|N[A]|}
  = \frac{\alpha(G\times H)}{|V(G\times H)|}=\frac{\alpha(G)}{|V(G)|}$.
If $E(G)=\emptyset$, the result is trivial, so we suppose
$E(G)\neq\emptyset$, then Corollary  \ref{bi} implies that
  $\frac{\alpha(H)}{|V(H)|}\leq \frac{\alpha(G)}{|V(G)|}\leq \frac{1}{2}$.
  By Proposition \ref{coro1}, there exists some $S\in
  I(G\times H)$ such that $A=S\cap N[A]$. Since $G\times H$ is MIS-normal,  we may
   assume that $S=S'\times H$ for some $S'\in I(G)$. Thus $A=(S'\times H)\cap N[A]$.
Set $B=\overline{N}[A]$. Then, by Proposition \ref{irregular},
$\sigma(B)\cap B=\emptyset$ or $B$ for every $\sigma\in
\mbox{Aut}(G\times H)$.

Set $C=\partial_G(B)$. For every pair  $a$ and $b$ of $C$,  select
$u\in
\partial_H(a,B)$ and $v\in
\partial_H(b,B)$. Since  $G$ and $H$ are vertex-transitive, there
exist  $\gamma\in \mbox{Aut}(G)$ and $\tau\in \mbox{Aut}(H)$ such
that $a=\gamma(b)$ and $u=\tau(v)$. It is clear that
$\sigma=(\gamma,\tau)\in \mbox{Aut}(G\times H)$ and
$(a,u)=\sigma(b,v)\in \sigma(B)\cap B$. By Proposition
\ref{irregular}, we conclude that $\sigma(B)=B$. Thus, we have
$\partial_H(a,B)=\tau(\partial_H(b,B))$. Therefore,
$|\partial_H(a,B)|=|\partial_H(b,B))|$  for any $a, b\in C$. In the
following, we will complete the proof by two cases.

\textit{Case $1$}: $C\neq V(G)$. Set $\overline{C}=(V(G)-C)$. Then
$(\overline{C}\times H)\cap B=\emptyset$, and thus
$\overline{C}\times H\subseteq N[A]$. For every $S''\in I(G)$, it is
clear that $S''\times H$ is a maximum independent set of $G\times
H$. Since $\frac{\alpha (B)}{|B|}=\frac{\alpha(G\times H)}{|G\times
H|}=\frac{\alpha(G)}{|V(G)|}$ and
$|\partial_H(a,B)|=|\partial_H(b,B)|$ for all $a,b\in
\partial_G(B)$, from Lemma \ref{caku} and the MIS-normality of $G\times H$
it follows that
$$\frac{|(S''\times H)\cap B|}{|B|}=\frac{|S''\cap
C|}{|C|}=\frac{\alpha(G)}{|V(G)|}.$$ Thus for every  $S''\in I(G)$,
\begin{eqnarray}\label{eq1}\frac{\alpha(G)}{|V(G)|}=\frac{|S''|}{|V(G)|}=\frac{|S''\cap
C|+|S''\cap \overline{C}|}{|C|+|\overline{C}|}=\frac{|S''\cap
\overline{C}|}{|\overline{C}|}=\frac{|S''\cap
C|}{|C|}.\end{eqnarray} Recall that $\overline{C}\times H\subseteq
N[A]$ and  $A\subseteq S'\times H$, it is easy to see  that
$A=N[A]\cap (S'\times H)$  and $\partial_G(A\cap (\overline{C}\times
H))=S'\cap \overline{C}$. Setting
 $F=S'\cap \overline{C}$, we have that $a\times H\subseteq A$ for every
$a\in F$. If $N_G[F]\cap C\neq\emptyset$, then there exist $a\in F$
and $b\in {C}$ such that $(a,b)\in E(G)$.  Since $B=\overline{N}[A]$
and $a\times H\subseteq A$, by definition,
$(b,u)\subseteq\overline{N}[a\times H]$ for every $u\in
\partial_H(b,B)$. Hence $ N_H[H]\neq \emptyset$ and $E(H)=\emptyset$, which
contradicts that $\frac{\alpha(H)}{|H|}\leq \frac{1}{2}$. Thus
$N_G[F]\cap C=\emptyset$, i.e., $N_G[F]\subseteq \overline{C}$. By
Proposition \ref{coro1}  and  (\ref{eq1}),
$$\frac{\alpha(G)}{|V(G)|}\geq \frac{|F|}{|N_G[F]|}=\frac{|
{S'}\cap \overline{C}|}{|N_G[F]|}\geq \frac{|{S'}\cap
\overline{C}|}{|\overline{C}|}=\frac{\alpha(G)}{|V(G)|}.$$ Therefore
$\frac{|F|}{|N_G[F]|}=\frac{\alpha(G)}{|V(G)|}$, so $G$ is
IS-imprimitive and (i) holds.

\textit{Case $2$}: $C=V(G)$. Since
$|\partial_H(a,B)|=|\partial_H(b,B))|$  for all $a, b\in V(G)$, we
have  $\partial_G(N[A])=V(G)$ and $|\partial_H(a,
N[A])|=|\partial_H(b, N[A])|< |H|$ for all $a,b\in V(G)$. Since
$A=(S'\times H)\cap N[A]$,  $\partial_H(a,N[A])\subseteq
\partial_H(a,S'\times H)$ for all $a\in \partial_G(A)$. Thus $\partial_H(a,A)=\partial_H(a, N[A])$ for all
$a\in
\partial_G(A)$. Select two
vertices $a$ and $b$ of $V(G)$ such that $a\in
\partial_G(A)$ and $(a,b)\in E(G)$. Then, for every $u\in [V(H)-\partial_H(b,N[A])]$ and $v\in \partial_H(a,N[A])$,
it is clear that $[(b, u), (a,v)]\not\in E(G\times H)$, so
$(u,v)\not\in E(H)$. This means $u\not\in N_H(\partial_H(a,N[A]))$,
that is,
\begin{eqnarray}\label{l}V(H)-\partial_H(b,N[A])\subseteq V(H)- N_H(\partial_H(a,N[A])).\end{eqnarray}
 If $\partial_H(b,N[A])=\partial_H(a,N[A])$, it follows from (\ref{l}) that $H$ is disconnected, and so either (i) or (ii) holds.

 Suppose that $\partial_H(b,N[A])\neq\partial_H(a,N[A])$ and  set
$D=\partial_H(a,N[A])-\partial_H(b,N[A])$.
  It is easy to check that
\begin{eqnarray*}2|D|=
|\partial_H(a,N[A])\cup\partial_H(b,N[A])-
\partial_H(a,N[A])\cap\partial_H(b,N[A])|.
\end{eqnarray*}
Since $D\subseteq H-\partial_H(b,N[A])$ and $D\subseteq
\partial_G(a,N[A])$, by (\ref{l}), we have $$D\subseteq V(H)-\partial_H(b,N[A])\subseteq V(H)-{N}_H(\partial_H(a,N[A]))\subseteq
V(H)-{N}_H(D).$$ So $D$ is an independent set of $H$ and
\begin{eqnarray*}N_{ H}[D]&\subseteq& D\cup
[\partial_H(b,N[A])-\partial_H(a,N[A])]\\
&=&\partial_H(a,N[A])\cup\partial_H(b,N[A])-\partial_H(a,N[A])
\cap\partial_H(b,N[A]),\end{eqnarray*} which implies that
$\frac{1}{2}\geq \frac{\alpha(H)}{|V(H)|}\geq \frac{|D|}{|N_{
H}[D]|}\geq\frac{1}{2}$. Thus
$\frac{\alpha(G)}{|V(G)|}=\frac{\alpha(H)}{|V(H)|}=\frac{1}{2}$. By
Corollary \ref{bi}, $G$ and $H$ are both bipartite, so (i) holds and the proof  completed. \qed

\begin{thm}\label{111}
Let $G$ and $H$ be two non-bipartite vertex-transitive graph such
that $\frac{\alpha(H)}{|V(H)|}=\frac{\alpha(G)}{|V(G)|}$. If
$G\times H$ is MIS-normal, then  $G$, $H$ and $G\times H$ are all
IS-primitive.
\end{thm}
\textbf{Proof.} First, suppose that  $G$  is IS-imprimitive and let
$A$ be an imprimitive independent set in $G$.
 For any $S\in I(H)$, let $S'=(\overline{N}_{G}[A]\times S)\cup (A\times H)$. It is clear that  $S'$ is an independent set of $G\times H$ and
\begin{eqnarray*}
|S'|&=&|\overline{N}_{G}[A]\alpha(H)|+|A||V(H)|=(|\overline{N}_{G}[A]|+|N_{G}[A]|)\alpha(H)\\
&=&|V(G)|\alpha(H)=\alpha(G\times H),
\end{eqnarray*}
i.e., $S'$ is a maximum independent set of $G\times H$,
contradicting the MIS-normality of $G$. Therefore, $G$ is
IS-primitive. Similarly, $H$ is also IS-primitive. By Lemma
\ref{dgp}, $G\times H$ is IS-primitive.\qed

\section{MIS-normality of the Products of Graphs}

 The following theorem is the main result on the MIS-normality of
 products of vertex-transitive graphs in this paper.
\begin{thm}\label{irregular product}
Let  $G$ and $H$ be two vertex-transitive graphs. Suppose that there
exists an induced subgraph $G'$ of $G$ such that
 $G'\times H$ is MIS-normal and $\frac{\alpha(G')}{|V(G')|}=\frac{\alpha(G)}{|V(G)|}$.
 Then either: (i) $G\times H$ is  MIS-normal, or (ii) $\frac{\alpha(G)}{|V(G)|}=\frac{\alpha(H)}{|V(H)|}$ and
 $G$ is IS-imprimitive, or (iii)  $\frac{\alpha(G)}{|V(G)|}<\frac{\alpha(H)}{|V(H)|}$ and $G$ is
disconnected.
\end{thm}
\textbf{Proof.} If $E(H)=\emptyset$, the result is obvious, so we
assume that $E(H)\neq\emptyset$. By Lemma \ref{caku} and the
MIS-normality of $G'\times H$, we have the following inequality
$$\frac{\alpha(G\times H)}{|V(G)||V(H)|}\leq \frac{\alpha(G'\times
H)}{|V(G')||V(H)|}=\max\left\{\frac{\alpha(G)}{|V(G)|},
\frac{\alpha(H)}{|V(H)|}\right\}\leq \frac{\alpha(G\times
H)}{|V(G)||V(H)|},$$ yielding
\begin{eqnarray}\label{33}\frac{\alpha(G\times
H)}{|V(G)||V(H)|}=\frac{\alpha(G'\times
H)}{|V(G')||V(H)|}=\max\left\{\frac{\alpha(G)}{|V(G)|},
\frac{\alpha(H)}{|V(H)|}\right\}.\end{eqnarray}

For every $\sigma\in \mbox{Aut}(G)$, it is clear that
  $\sigma(G')\times H$ is  MIS-normal.
Let $S$ be a maximum independent set of $G\times H$. By Lemma
\ref{caku} and (\ref{33}),  $S\cap (\sigma(G')\times H)$ is a
maximum independent set of $\sigma(G')\times H$. Clearly, for each
$a\in
\partial_G(S)$, there is a $\sigma\in \mbox{Aut}(G)$ such that $a\in
\sigma(G')$. We therefore have that $|\partial_H(a,S)|=|H|$ or
$\alpha(H)$ for each $a\in
\partial_G(S)$. In the following we distinguish three cases to
complete the proof.

\textit{ Case $1$}:  $|\partial_H(a,S)|=|V(H)|$ for every $a\in
\partial_G(S)$. By (\ref{33}),  we obtain that
$|\partial_G(S)|=\alpha(G)$. Since $E(H)\neq\emptyset$,
$\partial_G(S)$ is  an independent set of $G$. This implies that
$S=\partial_G(S)\times H$.

\textit{ Case $2$}: $|\partial_H(a,S)|=\alpha(H)$ for every $a\in
\partial_G(S)$. By (\ref{33}), we have that $\partial_G(S)=G$, $\frac{\alpha(H)}{|V(H)|}\geq \frac{\alpha(G)}{|V(G)|}$ and
$\partial_H(a,S)$ is a maximum independent set of $H$ for every
$a\in G$. Let $a$ be a fixed vertex of $G$, and set
$$C=\{c\in G:
\partial_H(c,S)=\partial_H(a,S)\}.$$
If $C=G$, then $S=G\times
\partial_H(a,S)$. If $C\neq G$, then choose $d\in G-C$ and $c\in C$. Since $\partial_H(c,
S)\neq \partial_H(d, S)$, there exists $u\in\partial_H(c, S)$ and $v
\in \partial_H(d, S)$ such that $(u, v)\in E(H)$ and $[(c,u),(d, v)]
\not\in E(G\times H)$. This implies that $(c, d) \in  E(G)$ and thus
G is disconnected.

\textit{ Case $3$}: $|\partial_H(a,S)|=|V(H)|$ and
$|\partial_H(b,S)|=\alpha(H)$ for some $a,b\in \partial_G(S)$. By
(\ref{33}), $\frac{\alpha(H)}{|V(H)|}= \frac{\alpha(G)}{|V(G)|}$ and
$\alpha(G\times H)= \alpha(G)|V(H)|=\alpha(H)|V(G)|$.  Set
$$C=\{c \in G:|\partial_H(c,S)|=|V(H)|\} \mbox{\ and \ }D=\{d\in G:
|\partial_H(d,S)|=\alpha(H)\}.$$ Since $E(H)\neq \emptyset$, it is
clear that $C$ is an independent set of $G$ and $(c,d)\not\in E(G)$
for every $c\in C$ and $d\in D$. So $N_G[C]\subseteq V(G)-D$.
Moreover,
$$|S|=\alpha(H)|V(G)|=|C||V(H)|+|D|\alpha(H).$$
Thus $\frac{|C|}{|N_G[C]|}\geq
\frac{|C|}{|V(G)|-|D|}=\frac{\alpha(H)}{|V(H)|}=
\frac{\alpha(G)}{|V(G)|}$. By Proposition \ref{coro1},
$\frac{|C|}{|N_G[C]|}=\frac{\alpha(G)}{|V(G)|}$, that is,   $G$ is
IS-imprimitive.

This completes the proof.\qed

The following Corollary solves Problem \ref{prb1} in a bit more
general setting.
\begin{corollary}\label{Mthm}
Let $G$ be a  vertex-transitive, non-bipartite graph. If $G^2$ is
MIS-normal, then $G^n$ is also MIS-normal and IS-primitive for all
$n\geq 3$.
\end{corollary}

\noindent\textbf{Proof}. We prove by induction on $n$. Since $G^2$
is MIS-normal, by Theorem \ref{111},
   $G$ and $G^2$
are both IS-primitive. Assume that $G^d$ is
 MIS-normal and IS-primitive for
all $d=2,\ldots,n-1$. We now prove that  $G^n$ is  MIS-normal and
IS-primitive. Note that $G^n=G^2\times G^{n-2}$. Let $G'$ be some
subgraph of $G^2$ that is isomorphic to $G$, for instance, the
subgraph induced by the set of vertices $\{(u, u) : u \in V (G)\}$.
It is clear that
$$\frac{\alpha(G')}{|V(G')|}=\frac{\alpha(G)}{|V(G)|}=\frac{\alpha(G^2)}{|V(G^2)|}$$
and  $G'\times G^{n-2}$ is isomorphic to $G^{n-1}$.
 Thus by assumption, $G'\times G^{n-2}$ is MIS-normal. By Theorem \ref{irregular
product} and Theorem \ref{111}, it is easy to see that $G^n$ is
MIS-normal and  IS-primitive. This completes the proof. \qed

\noindent\textbf{Acknowledgement} The author is greatly indebted to
the anonymous referees for giving
  useful comments and suggestions that have considerably improved the
  manuscript.
  He is grateful also for many valuable discussions with Professor J. Wang and Professor C.J. Zhou.


\begin{thebibliography}{99}

\bibitem{makl}
M.O. Albertson and K.L. Collins, Homomorphisms of $3$-chromatic
graphs, Discrete Math., 54 (1985) 127-132.

\bibitem{larose}
B. Larose and C. Tardif, Projectivity and independent sets in powers
of graphs, J. Graph Theory, 40 (2002) 162-171.


\bibitem{Tardif}
C. Tardif, Graph products and the chromatic difference sequence of
vertex-transitive graphs, Discrete Math., 185 (1998) 193-200.


\bibitem{Cameron}
P.J. Cameron and C.Y. Ku, Intersecting families of permutations,
European J. Comb., 24 (2003) 881-890.



\bibitem{MJ}
 V.P. Mario and J. Vera, Independent and coloring properties of direct
products of some vertex-transitive graphs, Discrete Math., 306
(2006) 2275-2281.


\bibitem{frankl3}
P. Frankl, An Erd\H{o}s-Ko-Rado Theorem for direct products,
European J. Combin., 17 (1996) 727-730.

\bibitem{Kla}
P.K. Jha and S. Klav\u{z}ar, Independence in direct-product graphs,
Ars Combin., 50 (1998) 53-60.


\bibitem{Ah1}
R. Ahlswede,  H. Aydinian and L.H. Khachatrian, The intersection
theorem for direct products, European J. Combin., 19 (1998) 649-661.







\end{thebibliography}
\end{document}